\documentclass{ws-ijnt}
\usepackage{amsmath}
\usepackage{amsfonts}

\newcommand\OO{{\mathcal O}}
\newcommand\NN{{\mathcal N}}
\newcommand\UU{{\mathcal U}}

\newcommand\bx{{\mathbf x}}
\newcommand\genus{{\mathbf g}}

\newcommand\A{{\mathbb A}}

\newcommand\PP{{\mathbb P}}
\newcommand\Q{{\mathbb Q}}

\newcommand\Z{{\mathbb Z}}

\makeatletter
\renewcommand{\p@enumii}{}
\makeatother

\title{CHARACTERIZING ALGEBRAIC CURVES WITH INFINITELY MANY INTEGRAL POINTS}

\author{PARASKEVAS ALVANOS}
\address{Aristotle University of Thessaloniki, Department of Mathematics\\
54124 Thessaloniki, Greece\\ 
paris14@math.auth.gr}

\author{YURI BILU}
\address{Universit\'e Bordeaux~1, Institut de Math\'ematiques\\
351 cours de la Lib\'eration, 33405 Talence, France\\
yuri@math.u-bordeaux1.fr}

\author{DIMITRIOS POULAKIS}
\address{Aristotle University of Thessaloniki, Department of Mathematics\\
54124 Thessaloniki, Greece\\ 
poulakis@math.auth.gr}

\begin{document}

\markboth{P. Alvanos, Yu. Bilu, D. Poulakis}
{Curves with infinitely many integral points}

\catchline{}{}{}{}{}

\renewcommand\thefootnote{\arabic{footnote}}

\maketitle

\begin{history}
\received{10.10.2007}
\accepted{+++++++++++++}
\comby{Umberto Zannier}
\end{history}

\begin{abstract}
\noindent
A classical theorem of Siegel asserts that the set of $S$-integral
points of an algebraic curve~$C$ over a number field  is finite
unless~$C$ has genus~$0$ and at most two points at infinity.
In this paper we give necessary and sufficient
conditions for~$C$ to have infinitely many $S$-integral points.
\end{abstract}


\section{Introduction}

Let ${C\subset \A^n}$ be an affine algebraic curve defined over a number field~$K$. By ${x_1, \ldots, x_n}$ we denote the coordinate functions on the affine space~$\A^n$ and their restrictions to~$C$, so that ${K(C)=K(x_1, \ldots, x_n)}$. 

By the \emph{points at infinity} of~$C$ we mean the infinite places of the field $\bar K(C)$ (that is, places where at least one of the coordinate functions ${x_1, \ldots, x_n}$ admits a pole). We denote the set of points at infinity by~$C_\infty$. By the \emph{genus} $\genus(C)$  we mean the genus of the field $\bar K(C)$ of $\bar K$-rational functions on~$C$ (or, equivalently, the genus of any non-singular projective model of~$C$). 

Let~$S$ be a finite set of places of~$K$ containing the set~$S_\infty$ of infinite places, and let ${\OO_S=\OO_{K,S}}$ be the set of $S$-integers of~$K$. We denote by $C(\OO_S)$ the set of points on~$C$ with coordinates in~$\OO_S$:
$$
C(\OO_S)= C(K) \cap \A^n(\OO_S). 
$$
The classical theorem of Siegel~\cite{Si29} asserts that the set $C(\OO_S)$ is finite if ${\genus(C)>0}$ or if~$C$ has at least~$3$ points at infinity (see also \cite[Chapter~8]{La83} and \cite[Chapter~D.9]{HS00}). However, this set can be finite even if ${\genus(C)=0}$ and ${|C_\infty|\le 2}$. The purpose of this note is to give a reasonable necessary and sufficient condition
for infinitude of the set $C(\OO_S)$. 

\begin{theorem}
\label{thos}
In the set-up above, the following two conditions are equivalent.
\begin{enumerate}
\item
\label{if}
The set $C(\OO_S)$ is infinite. 

\item
\label{ic}
The curve~$C$ is of genus~$0$, the set $C(\OO_S)$ contains a non-singular point and the set~$C_\infty$ has one of the following properties. 

\begin{enumerate}
\item
\label{i1}
We have ${|C_\infty|=1}$. 

\item 
\label{i2k}
We have ${|C_\infty|=2}$, both the points at infinity are defined over~$K$ and ${|S|\ge 2}$.

\item
\label{i2l}
We have ${|C_\infty|=2}$, the points at infinity are conjugate over~$K$, and at least one ${v\in S}$ splits in the field of definition of the points at infinity (which is quadratic over~$K$).
\end{enumerate}
\end{enumerate}
\end{theorem}

In the important special case, when ${S=S_\infty}$ is the set of infinite places and ${\OO_S=\OO}$ is the ring of integers of the field~$K$, we obtain the following statement.

\begin{theorem}
\label{tho}
The following two conditions are equivalent. 
\begin{enumerate}
\item
The set $C(\OO)$ is infinite. 

\item
The curve~$C$ is of genus~$0$, the set $C(\OO)$ contains a non-singular point and the set~$C_\infty$ has one of the following properties. 

\begin{enumerate}
\item
We have ${|C_\infty|=1}$. 

\item 
We have ${|C_\infty|=2}$, both points at infinity are defined over~$K$ and~$K$ is neither~$\Q$ nor an imaginary quadratic field.

\item
We have ${|C_\infty|=2}$, the points at infinity are conjugate over~$K$, and the field of definition of the points at infinity  is not a CM-extension of~$K$.  
\end{enumerate}
\end{enumerate}
\end{theorem}

Recall that a quadratic extension of number fields $L/K$ is called a \emph{CM-extension}  if the field~$K$ is totally real and~$L$ is  totally imaginary.

Theorem~\ref{tho} is an immediate consequence of Theorem~\ref{thos}. Indeed, ${|S_\infty|=1}$ if and only if~$K$ is either~$\Q$ or an imaginary quadratic field. Also, if~$L$ is a quadratic extension of~$K$ such that no infinite place of~$K$ splits in~$L$, then every infinite place of~$K$ is real and becomes complex in~$L$, which exactly means that $L/K$ is a CM-extension. 

When~$C$ is a line or a conic in~$\A^2$, Theorem~\ref{thos} was proved 
by Beukers~\cite{Be95}. For the case ${K=\Q}$ and ${\OO_S=\Z}$ see \cite[Theorem~A]{Si00}, \cite[Theorem~5.2]{PV02} and~\cite{Po03}. In~\cite{PV02} the problem is also approached from the computational viewpoint.

In~\cite{Si00}
Silverman writes: \textit{Little of the material in this note will be new to the ``experts'', but
we hope that an elementary exposition will be a useful addition to
the literature.} We believe that the same statement is applicable to the present article. While
our Theorem~\ref{thos}  will likely be unsurprising to the ``experts'', it does
not seem to have ever appeared before in this generality, and we believe that 
having it in the mathematical
literature will be quite useful.

A sort of generalization of Theorem~\ref{thos} to complements of hyperplanes on~$\PP^n$ has recently been obtained by Levin \cite[Theorem~3]{Le08}.

\section{Proof of Theorem~\ref{thos}}
We start from the \textbf{implication {\ref{if}$\Rightarrow$\ref{ic}}}. 
If $C(\OO_S)$ is infinite then the theorem of Siegel implies that ${\genus(C)=0}$ and ${|C_\infty|\le 2}$. Also, since there can be only finitely many singular points, there exist infinitely many non-singular points in $C(\OO_S)$. 

We have to show that one of the properties~(\ref{i1}),~(\ref{i2k}) or~(\ref{i2l}) is satisfied. If ${|C_\infty|=1}$ then we have~(\ref{i1}). Now assume that  ${|C_\infty|=2}$, and write ${C_\infty = \{A,B\}}$. Two cases are possible: either \emph{both~$A$ and~$B$ are defined over~$K$}, or \emph{they are conjugate over~$K$}. 

\medskip

Assume that \textbf{$A$ and~$B$ are defined over~$K$}. 
Since ${\genus(C)=0}$, there exists ${u\in K(C)}$ such that ${(u)=A-B}$ (where~$(u)$ is the principal divisor of~$u$). Both~$u$ and~$u^{-1}$ are integral over the ring ${K[\bx]=K[x_1, \ldots, x_n]}$. Hence there exists a positive integer~$N$ such that both $Nu$ and $Nu^{-1}$ are integral over $\OO_S[\bx]$. 

Now let~$P$ be a non-singular point from $C(\OO_S)$. Since it is non-singular, it corresponds to a single place of $K(C)$. Hence  the value $u(P)$ is well-defined, and both $Nu(P)$ and $Nu(P)^{-1}$ are integral over ${\OO_S[\bx(P)]=\OO_S}$. But the ring~$\OO_S$ is integrally closed, whence ${Nu(P), Nu(P)^{-1} \in \OO_S}$. It follows that we have only finitely many possibilities for the fractional ideal $(u(P))$ in the ring~$\OO_S$. In other words, we have ${u(P)=\alpha\eta}$, where~$\alpha$ belongs to a finite set, and~$\eta$ is an $S$-unit. But if ${|S|=1}$ then the group of $S$-units is finite, which leaves only finitely many possibilities for $u(P)$, a contradiction. Hence ${|S|\ge 2}$, and we we have Property~(\ref{i2k}). This completes the proof in the case when both~$A$ and~$B$ are defined over~$K$. 

\medskip

Finally, assume that~\textbf{$A$ and~$B$ are  conjugate over~$K$}. Then ${K(A)=K(B)=L}$ is a quadratic extension of~$K$. We denote by~$T$ the set of places of the field~$L$, extending the places from~$S$, and by ${\OO_T=\OO_{L,T}}$ the ring of $T$-integers of the field~$L$.  We select~$u$ as above, but now ${u\in L(C)}$. Multiplying~$u$ by a suitable integer, we may assume that it is integral over the ring $\OO_S[\bx]$. It follows that ${u(P)\in \OO_T}$ for every non-singular ${P\in C(\OO_T)}$. 

Let~$\tau$ be the non-trivial automorphism  of $L/K$. Then~$\tau$ permutes~$A$ and~$B$, which implies  ${u^\tau=\lambda u^{-1}}$
for some ${\lambda \in L^\ast}$. 

We want to show that we have Property~(\ref{i2l}). Thus, let us assume that no ${v\in S}$ splits in~$L$ and derive a contradiction. 
 
Let ${P\in C(K)}$ be a non-singular $K$-rational point. Fix ${v\in S}$. By the assumption, it has a unique extension to~$L$; denote it by~$v$ as well. Then ${|\alpha^\tau|_v=|\alpha|_v}$ for any ${\alpha \in L}$. We apply this with ${\alpha =u(P)}$. Since ${u^\tau=\lambda u^{-1}}$ and  ${P^\tau=P}$ (because~$P$ is defined over~$K$), we have ${u(P)^\tau=\lambda u(P)^{-1}}$, and
$$
\left|u(P)\right|_v=\left|u(P)^\tau\right|_v=|\lambda|_v \left|u(P)\right|_v^{-1}, 
$$
which implies that ${|u(P)|_v=\sqrt{|\lambda|_v}}$.

Thus, $|u(P)|_v$ is bounded for any  ${v\in T}$. Since $u(P)$ is a $T$-integer,    this leaves only finitely many possibilities for $u(P)$, a contradiction. This completes the proof of the implication \ref{if}$\Rightarrow$\ref{ic}. 

\bigskip

Now let us prove the \textbf{implication \ref{ic}$\Rightarrow$\ref{if}}. 
Thus,  assume that ${\genus(C)=0}$, that~$C(\OO_S)$ has a non-singular point (denoted by~$P$) and that~$C_\infty$ has one of the three properties above. Let us show that $C(\OO_S)$ is infinite.

\medskip

\textbf{Assume~(\ref{i1})}, and write ${C_\infty=\{A\}}$. Since ${\genus (C)=0}$, there exists ${u\in K(C)}$ with ${(u)=P-A}$. This~$u$ generates the field $K(C)$, and, since the coordinate functions ${x_1, \ldots, x_n}$ have no poles other than~$A$, they all can be expressed as polynomials in~$u$:
$$
x_i=F_i(u) \qquad (i=1, \ldots, n),
$$
where\footnote{We use the small letter~$u$ to denote a rational function, and the capital~$U$ for an independent variable.}  ${F_i(U)\in K[U]}$. Since ${u(P)=0}$, we have 
$$
F_i(0)=x_i(P)\in \OO_S \qquad (i=1, \ldots, n).
$$
Let~$N$ be a positive integer such that ${NF_i(U)\in \OO_S[U]}$ for ${i=1, \ldots, n}$. Then for any ${\theta \in \OO_S}$, satisfying ${\theta \equiv 0\mod N}$, we have the congruence
$$
NF_i(\theta)\equiv NF_i(0)\equiv 0\mod N 
$$
in the ring~$\OO_S$. Hence ${F_i(\theta)\in \OO_S}$, and  ${\left(F_1(\theta),  \ldots, F_n(\theta)\right)}$ is an $S$-integral point on~$C$. Since there are infinitely many choices of~$\theta$, the set $C(\OO_S)$ is infinite. 

\medskip

Now \textbf{assume~(\ref{i2k})}, and write ${C_\infty=\{A, B\}}$. Then there exists ${u\in K(C)}$ with ${(u)=A-B}$, and we normalize it to have ${u(P)=1}$. This~$u$ again generates the field $K(C)$, and, since the coordinate functions ${x_1, \ldots, x_n}$ have no poles other than~$A$ and~$B$, they all can be expressed as \emph{Laurent polynomials}  in~$u$:
\begin{equation}
\label{elaur}
x_i=F_i(u) \qquad (i=1, \ldots, n),
\end{equation}
where ${F_i(U)\in K[U,U^{-1}]}$. Since ${u(P)=1}$, we have 
$$
F_i(1)=x_i(P)\in \OO_S \qquad (i=1, \ldots, n).
$$
Let~$N$ be a positive integer such that ${NF_i(U)\in \OO_S[U, U^{-1}]}$ for ${i=1, \ldots, n}$. Let~$\UU_S$ be the group of $S$-units of the field~$K$. This group is infinite because ${|S|\ge 2}$. Hence the kernel of the $\mod N$ reduction   map ${\UU_S\to (\OO_S/N\OO_S)^\times}$ is infinite as well. In other words, there exist infinitely many $S$-units ${\theta \in \UU_S}$ such that ${\theta \equiv 1\mod N}$. 
For any such~$\theta$  we have
$$
NF_i(\theta)\equiv NF_i(1)\equiv 0\mod N.
$$
Hence ${F_i(\theta)\in \OO_S}$, and  ${\left(F_1(\theta),  \ldots, F_n(\theta)\right)}$ is an $S$-integral point on~$C$. We again proved that the set $C(\OO_S)$ is infinite. 

\medskip

Finally, in \textbf{case~(\ref{i2l})}, we put ${L=K(A)=K(B)}$ and define ${u\in L(C)}$ by ${(u)=A-B}$ and ${u(P)=1}$. We again have~(\ref{elaur}), but now the Laurent polynomials $F_i(U)$ have coefficients in~$L$.

Let~$\tau$ be the non-trivial automorphism of $L/K$. We claim that ${u^\tau=u^{-1}}$. Indeed, since~$\tau$ permutes~$A$ and~$B$, we have ${u^\tau=\lambda u^{-1}}$ with some ${\lambda \in L^\ast}$. In particular, ${u(P)^\tau=\lambda u(P)^{-1}}$. Since ${u(P)=1}$, we have ${\lambda=1}$.

Further, since ${x_i\in K(C)}$, we have ${x_i^\tau=x_i}$, whence
$$
F_i(u)=F_i(u)^\tau= F_i^\tau(u^\tau)=F_i^\tau(u^{-1}).
$$
This proves  the Laurent polynomials~$F_i$ are ``skew-symmetric'', that is,  ${F_i^\tau(U)=F_i(U^{-1})}$. 

Let~$T$ be the set of places of~$L$ extending the places from~$S$, and let~$\OO_T$ and~$\UU_T$ be the ring of $T$-integers and the  group of $T$-units of~$L$, respectively. Since at least one place from~$S$ splits in~$L$, we have ${|T|>|S|}$. Hence the  rank of~$\UU_T$ is strictly bigger than that of~$\UU_S$.  It follows that the norm map ${\NN_{L/K}:\UU_T\to \UU_S}$ has infinite kernel. We denote this kernel by~$\UU_T^0$. For any ${\theta \in \UU_T^0}$ we have ${\theta^\tau=\theta^{-1}}$. Hence, for such a~$\theta$ we have 
$$
F_i(\theta)^\tau=F_i^\tau(\theta^\tau)=F_i\left((\theta^{-1})^{-1}\right)=F_i(\theta),
$$
that is, ${F_i(\theta)\in K}$. 

Now we complete the proof as in the previous case. 
Let~$N$ be a positive integer such that ${NF_i(U)\in \OO_T[U, U^{-1}]}$ for ${i=1, \ldots, n}$. Since the group~$\UU_T^0$ is infinite, so is  the kernel of the $\mod N$ reduction   map ${\UU_T^0\to (\OO_T/N\OO_T)^\times}$; that is, infinitely many  ${\theta \in \UU_T^0}$ satisfy ${\theta \equiv 1\mod N}$. 
For any such~$\theta$  we have
 ${F_i(\theta)\in \OO_S}$, and  ${\left(F_1(\theta),  \ldots, F_n(\theta)\right)}$ is an integral point on~$C$. Hence the set $C(\OO_S)$ is infinite. This proves Theorem~\ref{thos}. 
 
\section*{Acknowledgments}

We thank the anonymous referee for useful comments and suggestions. 

\medskip
\noindent
Paraskevas Alvanos was supported by the State Scholarships Foundation of Greece. Yuri Bilu was supported by the European ALGANT Scholarship.


\begin{thebibliography}{99}

\bibitem{Be95}
F. Beukers, Ternary form equations, \textit{J. Number Th.}~\textbf{54} (1995)
113–-133.

\bibitem{HS00}
M.~Hindry, J.~H.~Silverman,
Diophantine Geometry: an Introduction
(Springer Verlag,  2000).


\bibitem{La83} 
S.~Lang, \textit{Fundamentals of Diophantine Geometry} (Springer Verlag,  1983).


\bibitem{Le08}
A. Levin, On the Zariski-density of integral points on a complement
of hyperplanes in~$\PP^n$, \textit{J. Number Th.}~\textbf{128} (2008) 96--104.

\bibitem{PV02} D.~Poulakis and  E.~Voskos, Solving genus zero
Diophantine equations with at most two infinite valuations, \textit{J.~Symbolic Comput.} \textbf{33}(4) (2002) 479--491.

\bibitem{Po03} 
D.~Poulakis,
Affine curves with infinitely many integral points, \textit{Proc. Amer.
Math. Soc.} \textbf{131}(5) (2003) 1357--1359.

\bibitem{Si29}
C.~L.~Siegel,
\"Uber einige Anwendungen Diophantischer Approximationen,
\textit{Abh. Preuss Akad. Wiss. Phys.-Math. Kl.},
1929,
Nr.~1;
\textit{Ges. Abh.}, Band~1, 209--266.

\bibitem{Si00}
J.~H.~Silverman,
On the distribution of integer points on curves
of genus zero,
\textit{Theor. Computer Sc}.~\textbf{235} (2000) 163--170. 



\end{thebibliography}
\end{document}